\documentclass[11pt]{article}
\newfont{\bb}{msbm10}

\usepackage{amsfonts}
\usepackage{amsmath,amssymb}
\usepackage{booktabs}
\usepackage{amsbsy}
\usepackage{comment}
\usepackage{graphicx,color}
\usepackage{caption}
\usepackage{subcaption}

\def\Div{{\rm div}\,}
\newcommand{\bu} {\mathbf{u}}
\newcommand{\R} {\mathbb{R}}

\newcommand{\bn} {\mathbf{n}}
\newcommand{\tbn} {\widehat{\mathbf{n}}}

\newcommand{\bv} {\mathbf{v}}

\newcommand{\tbu} {\tilde{\mathbf{u}}}
\newcommand{\bw} {\mathbf{w}}
\newcommand{\bF} {\mathbf{F}}
\newcommand{\bS} {\mathbf{S}}
\newcommand{\bP} {\mathbf{P}}

\newcommand{\bx} {\mathbf{x}}

\newcommand{\bg} {\mathbf{g}}
\newcommand{\bq} {\mathbf{q}}
\newcommand{\QQ} {\mathbb{Q}}
\newcommand{\VV} {\mathbb{V}}

\newcommand{\bI} {\mathbf{I}}
\newcommand{\bE} {\mathbf{E}}
\newcommand{\bD} {\mathbf{D}}
\newcommand{\dO} {{\partial\Omega}}

\newcommand{\bA} {{\bf A}}

\newcommand{\T} {{\cal T}}

\def\tbu{ {\tilde \bu}}

\newcommand{\bsigma}{\mbox{\boldmath$\sigma$\unboldmath}}
\newcommand{\bpsi}{\mbox{\boldmath$\psi$\unboldmath}}
\newcommand{\bphi}{\mbox{\boldmath$\phi$\unboldmath}}
\newcommand{\bxi}{\mbox{\boldmath$\xi$\unboldmath}}

\begin{document}

\thispagestyle{empty}
\date{}

\title{A monolithic fluid-porous structure interaction finite element method  \thanks{This work has been supported by the Russian Science Foundation grant 19-71-10094.}}

\author{Alexander Lozovskiy\thanks{Marchuk Institute of Numerical Mathematics RAS; {\tt alex.v.lozovskiy@gmail.com}}
\and
Maxim A. Olshanskii\thanks{Department of Mathematics, University of Houston; {\tt molshan@math.uh.edu}}
\and
Yuri V. Vassilevski\thanks{Marchuk Institute of Numerical Mathematics RAS and Sechenov University; {\tt yuri.vassilevski@gmail.com}} }

\maketitle

\markboth{}{}

\begin{abstract} The paper introduces a fully discrete quasi-Lagrangian finite element method for a monolithic formulation of a fluid-porous structure interaction problem. The method is second order in time and allows a standard $P_2-P_1$ (Taylor--Hood) finite element spaces for fluid problems in both fluid and porous domains.  The performance of the method is illustrated on a series of numerical experiments.

\end{abstract}

%\begin{keywords}
%fluid-structure interaction, semi-implicit scheme, monolithic approach, blood flow,  numerical stability, finite element method
%\end{keywords}

%\begin{AMS} 76M10, 65M12, 74F10,  76Z05.
%\end{AMS}

%%%%%%%%%%%%%%%%%%%%%%%%%%%%%%%%%%%%%%%%%%%%%%%%%%%%%%%%%%%%%%

%%%%%%%%%%%%%%%%%%%%%%%%%%%%%%%%%%%%%%%%%%%%%%%%%%%%%%%%%%%%%%%%
%%%%%%%%%%%%%%%%%%%%%%%%%%%%%%%%%%%%%%%%%%%%%%%%%%%%%%%%%%%%%%%%
\section{Introduction}\label{intro}
Blood flow in a vessel with permeable walls or  penetration of  oil  through a crack in a porous matrix can be seen
as the interaction of a freely flowing fluid with a fluid-saturated poroelastic structure.
A continuum mechanics description  of such fluid-poroelastic phenomena often leads to  coupled systems  of (Navier--)Stokes and
Biot equations~\cite{showalter2005poroelastic,koshiba2007multiphysics}.
Recently, there has been a growing interest in the numerical solution of the Stokes--Biot and  Navier--Stokes--Biot problems.
Several authors suggested solution strategies based on decomposition of the system into fluid and poroelastic loosely coupled problems
to allow for a computationally efficient time-stepping schemes~\cite{badia2009coupling,bukavc2015partitioning}.
For the reason of better stability, monolithic  methods for the  (Navier--)Stokes--Biot equations have become popular in the literature.
They differ in the form of equations and the numerical treatment of the coupling conditions on the interface
between a free flow domain and a domain occupied by the porous structure.
In \cite{ambartsumyan2018lagrange} the continuity of fluid fluxes on the interface is imposed weakly with the help of a Lagrange multipier and
in~\cite{wen2020strongly} an interior penalty discontinuous Galerkin method is applied to obtain a discrete coupled formulation.
The Nitsche approach is used for coupling fluid and poroelastic finite element formulations in ~\cite{bukac2015effects,ager2019nitsche}.
Combination of the Nitsche approach and unfitted finite elements~\cite{ager2019nitsche} adds extra flexibility to the numerical solution.

Many  publications on numerical methods for the fluid--poroelastic problem ignore inertia effect in the fluid and
formulate the free fluid problem as a Stokes system. One reason for such  simplification is the lack of the energy dissipation principle
for the Navier--Stokes--Biot problem with the common interface conditions, which hinders the analysis in this case.
This issue is well-known already for the  Navier--Stokes--Darcy (the Navier--Stokes--Biot problem with rigid structure),
where a local well-posedness of the system is currently known only under a smallness  assumption (even in 2D) and
the proof uses involved arguments that work in the absence of \textit{a priori} energy bound~\cite{badea2010numerical,girault2009dg}.
In the context of the Navier--Stokes--Darcy coupling the issue was addressed in~\cite{ccecsmeliouglu2008analysis,cesmelioglu2013time},
where  interface conditions were modified to ensure the thermodynamical consistency of the complete system.
In this report, we follow~\cite{ccecsmeliouglu2008analysis,cesmelioglu2013time} and employ the suggested correction
to the stress balance in the  Navier--Stokes--Biot to end up with a dissipative system and stable numerical method.

We consider the Navier--Stokes--Biot system with the Beavers--Joseph--Saffman interface condition and a modified stress interface  condition
and discuss its energy balance. For an ALE formulation of the problem we further introduce a monolithic finite element method.
Our finite element method features the formulation of all equations in the reference coordinates encoding all information on geometry deformation
in solution-dependent coefficients. This formulation allows a simple application of the method of lines for the time discretization.
In particular, the second order discretization in time is straightforward.
Such monolithic approach was proved to be efficient for  FSI problems with an impermeable elastic
structure~\cite{Hron2006,lozovskiy2015unconditionally,lozovskiy2019analysis}, and we extend it here to the case of poroelasticity.
In the spirit of monolithic formulations we apply here the same finite elements to approximate fluid velocity and pressure in both domains.
We choose the Taylor--Hood element (P2-P1) for this purpose, which is a valid Darcy element for applications
where the local mass conservation is not critical~\cite{karper2009unified}. We use the same P2 element for the structure velocity.
To enforce the continuity of fluid flux through the interface, we use the penalty approach
(the  Nitsche approach as in~\cite{ager2019nitsche} would be an alternative).

The remainder of the paper is organized in three sections. We formulate the governing equations,
interface and boundary conditions in section~\ref{secModel}. The same section presents the integral formulation, the energy balance of the system,
and an ALE formulation that we use for the discretization. The finite element method is introduced in section~\ref{sec_FE}.
Section~\ref{s_num} presents results of several numerical experiments.

\section{FPSI model}\label{secModel}
Consider a time-dependent domain $\Omega(t)\subset\mathbb{R}^3$  containing fluid and an elastic porous structure.
A subdomain $\Omega^f(t)$  is entirely  occupied  by  fluid and a subdomian $\Omega^s(t)$ is occupied by porous elastic solid fully saturated with fluid.
These subdomains are non-overlapping and $\overline{\Omega(t)}=\overline{\Omega^f(t)}\cup\overline{\Omega^s(t)}$.
Two regions are separated by the interface $\Gamma^{fs}(t):=\dO^f(t)\cap\dO^s(t)$.

In this paper, the equations governing the fluid and solid motion will be written in the reference domains
\[
\Omega_f=\Omega^f(0),\quad \Omega_s=\Omega^s(0), \quad \Gamma_{fs}=\Gamma^{fs}(0).
\]
The deformation of the poroelastic part is given by the mapping
\[
\bxi_s~:~\Omega_s\times[0,t]\to\bigcup_{t\in[0,T]}\Omega^s(t),
\]
with the corresponding displacement $\bu_s$, $\bu_s(\bx,t):=\bx-\bxi_s(\bx,t)$ and \textit{ the velocity of the elastic structure}  $\bv_{s}=\partial_t\bu_s=-\partial_t\bxi_s(\bx,t)$.
\begin{figure}
  \centering
  \includegraphics[width=1\textwidth, trim=0 170 0 0, clip]{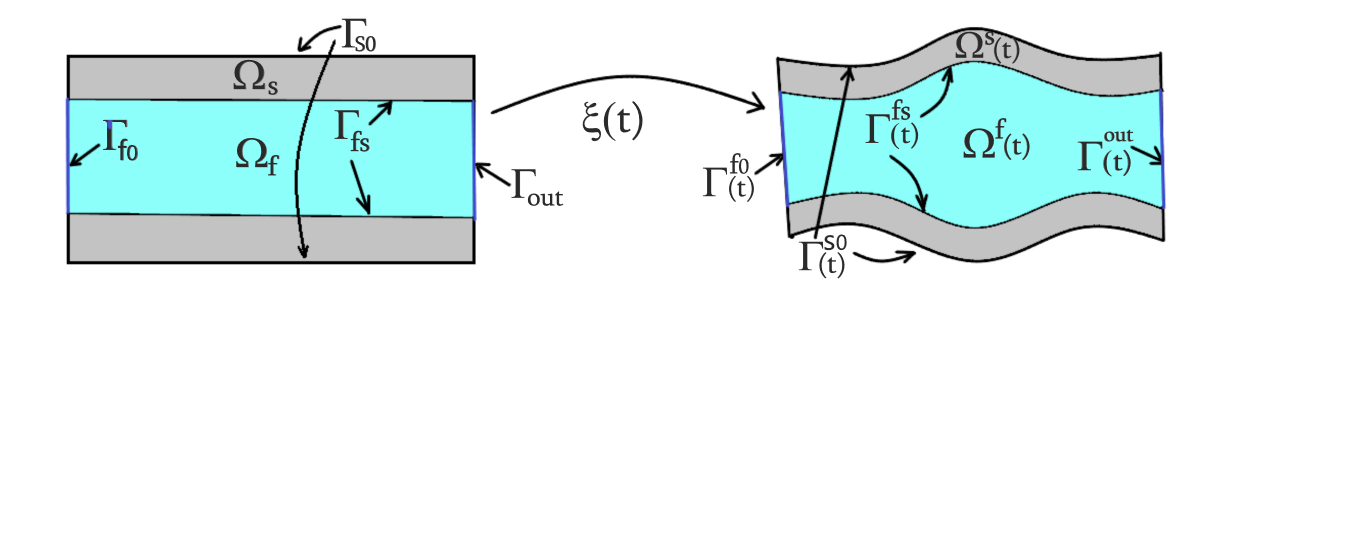}
  \caption{Reference and physical domains and boundaries.}\label{Fig1}
\end{figure}

The fluid dynamics is described by the velocity vector field $\bv(\bx,t)$ and the pressure function $p(\bx,t)$ defined in
the whole volume $\Omega(t)$ for all $t\in[0,T]$.
Following~\cite{koshiba2007multiphysics,badia2009coupling} we represent $\bv$ in the poroelastic domain
through the velocity of structure and the filtration flux $\bq=\phi(\bv-\bv_s)$, where $\phi$ is  the known porosity coefficient.
We denote the fluid pressure in the poroelastic domain by $p_d$, to emphasize its impact on the Darcy filtration, and
in the fluid domain by $p_f$.

Denote by $\rho_s$ and $\rho_f$ the densities of  solid and fluid. %Then $J\rho_s$ is the current density of the elastic structure.
Then $\rho_p=\rho_s(1-\phi)+\rho_f\phi$ is the density of the saturated porous medium.  Denote by $\bsigma_s$, $\bsigma_f$ the Cauchy stress
tensors in porous media and fluid, respectively. The poroelastic stress tensor is given by $\bsigma_p=\bsigma_s-\alpha p\bI$,
where $\alpha>0$ is  Biot's coefficient (typically $\alpha\simeq 1$, so further we set $\alpha=1$).
The porous medium is also characterized by its permeability tensor $K$.
The Biot system in the porous domain and the Navier-Stokes equations in the fluid domain follow from
the momenta balances and mass conservation principles (and neglecting the inertial effect of the matrix):

\begin{equation}\label{FPSI}
\begin{split}
&\left\{
\begin{split}
\rho_p\dot{\bv}_s+\rho_f\dot{\bq}&=\Div\bsigma_p\\
\rho_f\dot{\bv}_s+\frac{\rho_f}{\phi}\dot{\bq}&=-(K^{-1}\bq+\nabla p_d) \\
-s_0\dot{p_d}&=  \Div (\bv_s+\bq)
\end{split}\right. \quad\text{in}~\Omega^s(t),
\\
&\left\{
\begin{split}
\rho_f\dot{\bv}_f &=\Div \bsigma_f \\
\Div \bold{v}_f&= 0
\end{split}\right.\quad\text{in}~\Omega^f(t),
\end{split}
\end{equation}
where $1/s_0$ is Biot modulus or mixture compressibility modulus.

We divide the boundary of $\Omega(t)$  into the external boundary of the poroplastic   structure  $\Gamma^{s0}(t):=\dO(t)\cap\dO^s(t)$,
fluid Dirichlet and outflow boundaries: $\dO(t)\cap\dO^f(t)=\Gamma^{f0}(t)\cup\Gamma^{\rm out}(t)$; cf. Figure~\ref{Fig1}.
The governing equations are complemented with boundary conditions
\begin{equation}\label{FPSIbc}
\bold{v}^f = \bold{g}~ \text{ on }\Gamma^{f0}(t),\quad\bsigma_f\bn= \bold{0}\; \text{ on }\Gamma^{\rm out}(t),\quad p_d=0 \; \text{ on }\Gamma^{s0}(t),\quad \bold{v}^s = \bold{0}~ \text{ on }\Gamma^{s0}(t)
\end{equation}
and suitable initial conditions.

We now discuss coupling conditions on the interface between the fluid and poroelastic domains.
Denote by $\bn$ the normal vector on $\Gamma^{fs}(t)$ pointing from the fluid  to the poroelastic structure.
The balance of normal stresses on $\Gamma^{fs}(t)$  is commonly written in terms of  the  interface conditions:
$\bsigma_f\bn=\bsigma_p\bn$ and $\bn^T\bsigma_f\bn=-p_d$. This coupling, however, is not known to provide an energy consistent  (dissipative)
system. For the pure Darcy--Navier--Stokes coupling a remedy was suggested in \cite{ccecsmeliouglu2008analysis,ccecsmeliouglu2009primal}
where the second condition was changed to include a contribution of the fluid kinetic energy.
In this paper, we use the same modification in the poroelasticity context and  the two interface conditions read:
\begin{equation}\label{FPSIinterf}
\bsigma_f\bn%{\color{blue}+\frac{\rho_f}2(q\cdot\bn)_{-}\bv_f}
=\bsigma_p\bn\qquad\text{and}\quad \bn^T\bsigma_f\bn=-p_d +\frac{\rho_f}2|{\bv}_f|^2 ~~~\text{on}~\Gamma^{fs}(t).
\end{equation}
%The additional term on the left-hand side can be interpret as a correction to stress produced by the transport of momentum into the fluid region.
Such modification of the stress balance is similar to modifications of outflow boundary conditions
and 1D-3D models coupling conditions in computational fluid dynamics,  see e.g.,~\cite{bazilevs2009patient,braack2014directional}.
The continuity of the normal flux on the fluid-structure interface gives
\begin{equation}\label{cont_int}
\bv_f\cdot\bn=(\bv_s+\bq)\cdot\bn\quad\mbox{on}~\Gamma^{fs}(t).
\end{equation}
Finally, the Beavers--Joseph--Saffman condition sets the tangential component of the normal stress proportional
to the fluid ``slip''  rate along the interface:
\begin{equation}\label{BJSa}
\bP\bsigma_f\bn=-\gamma \bP K^{-\frac12}(\bv_f-\bv_s) ~~~\text{on}~\Gamma^{fs}(t),
\end{equation}
where $\bP$ is the orthogonal projector on the tangential plane to $\Gamma^{fs}(t)$.

\subsection{Integral formulation} In the preparation for the finite element method, we write out an integral (weak) formulation of the  FPSI problem \eqref{FPSI}--\eqref{BJSa}. We take the inner product of the elasticity equation in  \eqref{FPSI} with a sufficiently smooth $\bpsi_s$ such that $\bpsi_s=0$ on $\Gamma^{s0}(t)$, integrate it over $\Omega^s(t)$ and integrate the stress term by parts (recall that $\bn$ is inward for $\Omega^s$). This adds up with the first Darcy equation multiplied by a  sufficiently smooth $\bpsi_d$ and integrated over $\Omega^s(t)$ to give
\begin{multline}\label{aux194}
\int_{\Omega^s(t)} (\rho_p\dot{\bv}_s+\rho_f\dot{\bq})\cdot\bpsi_s + (\rho_f\dot{\bv}_s+\frac{\rho_f}{\phi}\dot{\bq}+ K^{-1}\bq)\cdot \bpsi_d\,dx
+ \int_{\Omega^s(t)}\bsigma_p:\nabla\bpsi_s\,dx \\
 - \int_{\Omega^s(t)}p_d\Div\bpsi_d     + \int_{\Gamma^{fs}(t)}(\bsigma_p\bn)\cdot\bpsi_s\,ds  -  \int_{\Gamma^{fs}(t)}p_d (\bpsi_d\cdot\bn)\,ds=0.
\end{multline}
Further, the fluid  momentum equation in \eqref{FPSI} is multiplied by a smooth vector function $\bpsi_f$ such that $\bpsi_f=\bold{0}$ on $\Gamma^{f0}$. Integrating over $\Omega^s(t)$ and integrating the stress term by parts we obtain
\begin{equation}\label{aux200}
\int_{\Omega^f(t)} \rho_f\dot{\bv}_f\cdot\bpsi_f\,dx+ \int_{\Omega^f(t)}\bsigma_f:\nabla\bpsi_f\,dx- \int_{\Gamma^{fs}(t)}\bpsi_f^T\bsigma_f\bn \,ds=0. %{\color{blue}+\frac{\rho_f}2(q\cdot\bn)_{-}\bv_f})\cdot\bpsi_f\,ds=0.
\end{equation}
We add up boundary terms in \eqref{aux194} and \eqref{aux200} and use interface conditions \eqref{FPSIinterf}--\eqref{BJSa} to reorganize them
\begin{align*}
  \int_{\Gamma^{fs}(t)}(\bsigma_p\bn)\cdot&\bpsi_s\,ds  -  \int_{\Gamma^{fs}(t)}p_d (\bpsi_d\cdot\bn)\,ds - \int_{\Gamma^{fs}(t)}\bpsi_f^T\bsigma_f\bn\,ds \\ %{\color{blue}+\frac{\rho_f}2(q\cdot\bn)_{-}\bv_f)\cdot\bpsi_f\,ds } \\
\text{\footnotesize use \eqref{FPSIinterf}} & = \int_{\Gamma^{fs}(t)}(\bsigma_f\bn)\cdot (\bpsi_s-\bpsi_f)\,ds -  \int_{\Gamma^{fs}(t)}p_d (\bpsi_d\cdot\bn)\,ds\\ % - {\color{blue}\int_{\Gamma^{fs}(t)}\frac{\rho_f}2(q\cdot\bn)_{-}\bv_f\cdot\bpsi_f\,ds}\\
 \text{\footnotesize split }\bsigma_f\bn &= \int_{\Gamma^{fs}(t)}(\bn^T\bsigma_f\bn) (\bpsi_s-\bpsi_f)\cdot\bn\,ds +\int_{\Gamma^{fs}(t)}(\bP\bsigma_f\bn)\cdot \bP(\bpsi_s-\bpsi_f)\,ds\\
 \qquad &-  \int_{\Gamma^{fs}(t)}p_d (\bpsi_d\cdot\bn)\,ds \\ %{\color{blue}- \int_{\Gamma^{fs}(t)}\frac{\rho_f}2(q\cdot\bn)_{-}\bv_f\cdot\bpsi_f\,ds}\\
  \text{\footnotesize use \eqref{FPSIinterf}, \eqref{BJSa}} &= \int_{\Gamma^{fs}(t)}p_d(\bpsi_f-\bpsi_s-\bpsi_d)\cdot\bn\,ds
  +\gamma\int_{\Gamma^{fs}(t)}K^{-\frac12}\bP(\bv_f-\bv_s)\cdot (\bpsi_f-\bpsi_s)\,ds\\
  &\qquad + \int_{\Gamma^{fs}(t)}\frac{\rho_f}2|{\bv}_f|^2  (\bpsi_s-\bpsi_f)\cdot\bn \,ds.
  %\\ &\qquad {\color{blue}- \int_{\Gamma^{fs}(t)}\frac{\rho_f}2(q\cdot\bn)_{-}\bv_f\cdot\bpsi_f\,ds}\\
\end{align*}
Summing up \eqref{aux194} and \eqref{aux200} and using the calculations above we arrive at the integral equality satisfied by
sufficiently smooth FPSI solution $\bv_s$, $\bq$, $\bv_f$, $p_d$, $p_f$
\begin{multline}\label{aux224}
	\int_{\Omega^s(t)} \left[ (\rho_p\dot{\bv}_s+\rho_f\dot{\bq})\cdot\bpsi_s + (\rho_f\dot{\bv}_s+\frac{\rho_f}{\phi}\dot{\bq}+ K^{-1}\bq)\cdot \bpsi_d
	\right]\,dx
+ \int_{\Omega^s(t)}\bsigma_p:\nabla\bpsi_s\,dx \\
 - \int_{\Omega^s(t)}p_d\Div\bpsi_d\,dx + \int_{\Omega^f(t)} \rho_f\dot{\bv}_f\cdot\bpsi_f\,dx+ \int_{\Omega^f(t)}\bsigma_f:\nabla\bpsi_f\,dx
% {\color{blue}- \int_{\Gamma^{fs}(t)}\frac{\rho_f}2(q\cdot\bn)_{-}\bv_f\cdot\bpsi_f\,ds}
  + \int_{\Gamma^{fs}(t)}\frac{\rho_f}2|{\bv}_f|^2  (\bpsi_s-\bpsi_f)\cdot\bn \,ds\\
   +\int_{\Gamma^{fs}(t)}p_d(\bpsi_f-\bpsi_s-\bpsi_d)\cdot\bn\,ds
  +\gamma\int_{\Gamma^{fs}(t)}K^{-\frac12}\bP(\bv_f-\bv_s)\cdot (\bpsi_f-\bpsi_s)\,ds =0
\end{multline}
for all sufficiently smooth $\bpsi_s$, $\bpsi_d$, and $\bpsi_f$ such that $\bpsi_s=0$ on $\Gamma^{s0}$, $\bpsi_f=\bold{0}$ on $\Gamma^{f0}$.
For the weak formulation, this integral identity should be supplemented by the two continuity equations in \eqref{FPSI}
and the normal continuity interface condition \eqref{cont_int}.

To obtain the energy balance identity,  we  assume that $\Gamma^{f0}$ and $\Gamma^{\rm out}$ are steady
and $\bg=0$ on $\Gamma^{f0}$.
We further let   $\bpsi_s=\bv_s$, $\bpsi_d=\bq$, $\bpsi_f=\bv_f$ and use $\bsigma_p=\bsigma_s- p\bI$,
continuity conditions and \eqref{cont_int} to arrive at the equality:
\begin{multline}\label{aux229}
	\int_{\Omega^s(t)} \left[(\rho_p\dot{\bv}_s+\rho_f\dot{\bq})\cdot\bv_s + (\rho_f\dot{\bv}_s+\frac{\rho_f}{\phi}\dot{\bq})\cdot\bq+ K^{-1}|\bq|^2\right]\,dx
+ \int_{\Omega^s(t)}\bsigma_s:\nabla\bpsi_s\,dx+\int_{\Omega^s(t)}s_0 \dot{p}_d p_d\,dx \\
    + \int_{\Omega^f(t)} \rho_f\dot{\bv}_f\cdot\bv_f\,dx+ \int_{\Omega^f(t)}\bsigma_f:\nabla\bv_f\,dx
 % {\color{blue}- \int_{\Gamma^{fs}(t)}\frac{\rho_f}2(q\cdot\bn)_{-}|\bv_f|^2\,ds}
  - \int_{\Gamma^{fs}(t)}\frac{\rho_f}2|{\bv}_f|^2  \bq\cdot\bn \,ds
  +\gamma\int_{\Gamma^{fs}(t)}K^{-\frac12}|\bP(\bv_f-\bv_s)|^2\,ds =0.
\end{multline}
Using $\bsigma_f=\mu_f\bD\bv_f-p_f\bI$, $\Div \bv_f=0$, and rearranging the first two term by substituting  $\rho_p=\rho_s(1-\phi)+\rho_f\phi$, we can rewrite the above equality as
\begin{multline}\label{aux258}
	\int_{\Omega^s(t)}\left[(1-\phi)\rho_s\dot{\bv}_s\cdot\bv_s +\phi\rho_f(\dot{\bv}_s+\frac{\dot{\bq}}{\phi})\cdot({\bv}_s+\frac{\bq}{\phi})\right]\,dx+ \int_{\Omega^s(t)}K^{-1}|\bq|^2\,dx
+ \int_{\Omega^s(t)}\bsigma_s:\nabla\bpsi_s\,dx\\
  +\int_{\Omega^s(t)}s_0 \dot{p}_d p_d\,dx  + \int_{\Omega^f(t)} \rho_f\dot{\bv}_f\cdot\bv_f\,dx+ \mu_f\int_{\Omega^f(t)}|\bD\bv_f|^2\,dx\\
%  {\color{blue}- \int_{\Gamma^{fs}(t)}\frac{\rho_f}2(q\cdot\bn)_{-}|\bv_f|^2\,ds}
- \int_{\Gamma^{fs}(t)}\frac{\rho_f}2|{\bv}_f|^2  \bq\cdot\bn \,ds
   +\gamma\int_{\Gamma^{fs}(t)}K^{-\frac12}|\bP(\bv_f-\bv_s)|^2\,ds =0.
\end{multline}
The integrals with material derivatives can be readily converted to the variations of kinetic energy by  application of the Reynolds transport theorem
and recalling that all parts of $\partial\Omega^f(t)$ are steady except $\Gamma^{fs}(t)$, which normal velocity is $\bv_s\cdot\bn$:
\[
\begin{split}
\frac{d}{dt}\frac12\int_{\Omega^f(t)} \rho_f|{\bv}_f|^2\,dx&=\int_{\Gamma^{fs}(t)} \rho_f\frac{\partial{\bv}_f}{\partial t}\cdot{\bv}_f\,ds
+ \frac12\int_{\Gamma^{fs}(t)} \rho_f|{\bv}_f|^2\bv_s\cdot\bn\,ds \\
&=\int_{\Gamma^{fs}(t)} \rho_f\frac{\partial{\bv}_f}{\partial t}\cdot{\bv}_f\,ds
+ \frac12\int_{\Gamma^{fs}(t)} \rho_f|{\bv}_f|^2\bv_f\cdot\bn\,ds - \frac12\int_{\Gamma^{fs}(t)} \rho_f|{\bv}_f|^2\bq\cdot\bn\,ds\\
&=\int_{\Gamma^{fs}(t)} \rho_f\frac{\partial{\bv}_f}{\partial t}\cdot{\bv}_f\,ds
+\frac12\int_{\Gamma^{fs}(t)} \rho_f\Div(|{\bv}_f|^2{\bv}_f)\,ds - \frac12\int_{\Gamma^{fs}(t)} \rho_f|{\bv}_f|^2\bq\cdot\bn\,ds\\
\text{\footnotesize using }  {\footnotesize\Div{\bv}_f=0}~ &=\int_{\Gamma^{fs}(t)} \rho_f\frac{\partial{\bv}_f}{\partial t}\cdot{\bv}_f\,ds
+\int_{\Gamma^{fs}(t)} \rho_f ((\bv_f\cdot\nabla){\bv}_f)\cdot{\bv}_f\,ds - \frac12\int_{\Gamma^{fs}(t)} \rho_f|{\bv}_f|^2\bq\cdot\bn\,ds\\
&=\int_{\Gamma^{fs}(t)} \rho_f\dot{\bv}_f\cdot{\bv}_f\,ds
 - \frac12\int_{\Gamma^{fs}(t)} \rho_f|{\bv}_f|^2\bq\cdot\bn\,ds.
\end{split}
\]
We handle the $\int_{\Omega^s(t)}$ integrals containing material derivatives in \eqref{aux258} by the same argument
assuming that the elastic structure is incompressible, i.e. $\Div\bv_s=0$, and recalling that the material derivative in
the structure is written in the Eulerian terms as $\partial/{\partial t} + \bv_s\cdot\nabla$.
%or additional terms accounting for the volume variation are included in to the equations of the Biot system in \eqref{FPSI}. In either case,
Therefore, \eqref{aux258} yields
\begin{multline}\label{aux235}
\frac{d}{dt}\frac12\int_{\Omega^s(t)}(1-\phi)\rho_s|\bv_s|^2 +\phi\rho_f|\bv_f|^2\,dx+ \int_{\Omega^s(t)}K^{-1}|\bq|^2\,dx
+ \int_{\Omega^s(t)}\bsigma_s:\nabla\bpsi_s\,dx\\
  +\frac{d}{dt}\frac12\int_{\Omega^s(t)}s_0 |{p}_d|^2\,dx  + \frac{d}{dt}\frac12\int_{\Omega^f(t)} \rho_f|{\bv}_f|^2\,dx+ \mu_f\int_{\Omega^f(t)}|\bD\bv_f|^2\,dx\\
  +\gamma\int_{\Gamma^{fs}(t)}K^{-\frac12}|\bP(\bv_f-\bv_s)|^2\,ds =0,
\end{multline}
where we  used $\bv_f=\bv_s+\frac{\bq}{\phi}$ in $\Omega^s(t)$ for the brevity.
We see that the system is  dissipative. % since the interface term $\int_{\Gamma^{fs}(t)}|{\bv}_f|^2(\bq\cdot\bn)_{+}\,ds$ is non-negative.
Without  the correction in the stress balance on the interface, the sign indefinite term
$-\frac{\rho_f}2\int_{\Gamma^{fs}(t)}|{\bv}_f|^2(\bq\cdot\bn)\,ds$ appears in the energy equality, and the system is not necessarily dissipative.
%A possible remedy is to change the second interface condition in \eqref{FPSIinterf} to include the kinetic energy:
%\begin{equation}\label{aux246}
% \bn^T\bsigma_f\bn=-p_d-\frac{\rho_f}2|{\bv}_f|^2 ~~~\text{on}~\Gamma^{fs}(t).
%\end{equation}
%In this case the extra term in \eqref{aux235} cancels out. Another option would be to update $\bsigma_f\bn$ in the interface conditions by including the influx of the energy: $\bsigma_f\bn+\rho_f(\bq\cdot\bn)_{-}\bu_f$.

\subsection{ALE formulation}

In this paper, we adopt the Arbitrary Lagrangian-Eulerian formulation by extending $\bxi_s$ to an  auxiliary mapping in the fluid domain
\[
\bxi_f~:~\Omega_f\times[0,t]\to\bigcup_{t\in[0,T]}\Omega^f(t)
\]
such that $\bxi_s=\bxi_f$ on $\Gamma_{fs}$, i.e. $\bxi$ is globally continuous. In general, $\bxi_f$ does not follow material trajectories. Instead, it is defined by
a continuous  extension of the displacement field to the flow reference domain
\begin{equation}\label{ext}
\bu_f:=\mbox{Ext}(\bu_s)=\bx-\bxi_f(\bx,t)\quad\mbox{in}~\Omega_f\times[0,t];\quad\bu = \left\{\begin{split}
                                                                 \bu_s & \text{ in } \Omega_s \\
                                                                 \bu_f & \text{ in } \Omega_f.
                                                              \end{split}\right.
\end{equation}
The corresponding globally defined deformation gradient is $\bF=\bI+\nabla\bu$, and $J:=\mbox{det}(\bF)$ is its determinant.
From now on, for notational simplicity, we will be using the same notation for these fields defined in the reference configuration as $\bv_f(\bx, t) := \bv_f(\bxi_f(\bx,t), t)$ and $p_f(\bx, t) := p_f(\bxi_f(\bx,t), t)$.
We use the notation  $\bsigma_s\circ\bxi_s(\bx):=\bsigma_s(\bxi_s(\bx))$.

The  governing equations driving the motion of fluid and  structure written in the reference domains read as
\begin{equation}\label{FSI1}
\left\{
\begin{split}
\rho_p\frac{\partial \bold{v}_s}{\partial t}+\rho_f\frac{\partial \bold{q}}{\partial t}=J^{-1}\Div(J(\bsigma_p\circ\bxi_s)\bF^{-T})\quad\text{in}~\Omega_s,\\
\rho_f\frac{\partial \bold{v}_s}{\partial t}+\frac{\rho_f}{\phi}\frac{\partial \bold{q}}{\partial t}=-K^{-1}\bq-\bF^{-T} \nabla p   \quad\text{in}~\Omega_s,\\
\rho_f\frac{\partial \bold{v}_f}{\partial t}=J^{-1}\Div(J(\bsigma_f\circ\bxi_f)\bF^{-T})- \rho_f\nabla \bold{v}_f (\bF^{-1}(\bold{v}_f-\frac{\partial \bold{u}}{\partial t}))
 \quad\text{in}~\Omega_f,
 \end{split}\right.
 \end{equation}
and the mass conservation reads as
\begin{equation}\label{FSI3}
\left\{
\begin{split}
\Div(J\bF^{-1}(\bold{v}_s+\bq))= -s_0J\frac{\partial p_d}{\partial t}&~~~\text{in}~\Omega_s,\\
\Div(J\bF^{-1}\bold{v}_f)= 0&~~~\text{in}~\Omega_f.
\end{split}
\right.
\end{equation}
Using the identity $\Div(J\bF^{-1}\bold{v})=J\nabla\bv:\bF^{-T}$, the last two equations can be written as
\begin{equation}\label{FSI3b}
\left\{
\begin{split}
\nabla(\bold{v}_s+\bq):\bF^{-T}= -s_0\frac{\partial p_d}{\partial t}&~~~\text{in}~\Omega_s,\\
\nabla\bv_f:\bF^{-T}= 0&~~~\text{in}~\Omega_f.
\end{split}
\right.
\end{equation}
The deformation of the structure can be found by integrating the kinematic equation
 \begin{equation}\label{FSI2}
\frac{\partial \bold{u}_s}{\partial t}=\bv_s~~~\text{in}~\Omega_s.
\end{equation}
The boundary and interface  conditions are the same in the ALE formulation.
The normal $\bn$ (and projector $\bP=\bI-\bn\bn^T$) to the interface and outflow boundary in the physical domain
can be computed from the reference normal $\tbn$, i.e. $\bn=\bF^{-T}\tbn/|\bF^{-T}\tbn|$. We collect all conditions in one place here:
\begin{equation}\label{FSIbc}
\bold{v}_f = \bold{g}~ \text{ on }\Gamma_{f0},\quad\bsigma_f\tbn= \bold{0}\; \text{ on }\Gamma_{\rm out},\quad \bold{v}_s = \bold{0}~ \text{ on }\Gamma_{s0},\quad p_d=0~ \text{ on }\Gamma_{s0}
%\bu=\bold{0}~ \text{ on }\Gamma_{s0}\cup\Gamma_{f0}\cup\Gamma_{\rm out},
\end{equation}
for the outer boundaries and
\begin{align}\label{FSIinterf}
\bsigma_f\bn %{\color{blue}+\frac{\rho_f}2(q\cdot\bn)_{-}\bv_f}
=\bsigma_p\bn,\qquad \bn^T\bsigma_f\bn=-p_d+\frac{\rho_f}2|{\bv}_f|^2 ~~~\text{on}~\Gamma_{fs},\\
\label{cont_int2}
\bv_f\cdot\bn=(\bv_s+\bq)\cdot\bn\quad\mbox{on}~\Gamma_{fs},\\
\label{BJS}
\bP\bsigma_f\bn=-\gamma \bP K^{-\frac12}(\bv_f-\bv_s) ~~~\text{on}~\Gamma_{fs}
\end{align}
on the interface. For the integral formulation in the reference coordinates,
we will  use the identities $J\,d\hat{x}=dx$, $J|\bF^{-T}\tbn|\,\,d\hat{s}=ds$,
where $ds$, $d\hat{s}$ are elementary areas orthogonal to $\bn$ and $\tbn$ in physical and reference coordinates, respectively.

The constitutive relation  for the Newtonian fluid in the reference domain reads
\begin{equation}\label{constit_f}
\bsigma_f=-p_f\bI+\mu_{f}(\nabla\bv \bF^{-1}+\bF^{-T}(\nabla\bv)^T)~~~\text{in}~\Omega_f.
\end{equation}
For the structure we consider  the compressible geometrically
nonlinear Saint Venant--Kirchhoff material
with
\begin{equation}\label{constit_s1}
\bsigma_s=\frac{1}{J}\bF\bS\bF^{T},\quad\text{with}~\bS=\lambda_s\mbox{tr}(\bE){\bI}+2\mu_{s}\bE,
\end{equation}
where  $\bE=\frac12\left(\bF^T\bF-\bI\right)$ is the  Lagrange-Green strain tensor and $\lambda_s,\mu_s$ are the Lame constants.

Thus, the FPSI problem in the reference coordinates consists in finding  pressure distributions $p_d$, $p_f$, fluid and structure  velocity fields $\bv_f$, $\bv_s$, fluid flux in the porous medium  $\bq$ and the displacement field
$\bu$   satisfying the set of equations, interface and boundary conditions
\eqref{FSI1}--\eqref{BJS}, together with  \eqref{constit_f}, \eqref{constit_s1}, and
subject to a given extension rule \eqref{ext}.
\medskip

\section{Discretization}\label{sec_FE}
We now proceed with  dicretization  of the FPSI problem formulated in the reference domain.
Treating the problem in the reference domain allows us to avoid time-dependent triangulations and finite element function spaces and
apply the standard method of lines to decouple space and time discretizations.
We adopt a finite element method in space  and define
 an admissible  triangulation of the reference domain $\overline{\Omega}(0)$ as a collection $\T_h$ of shape-regular tetrahedra
 such that the triangulation respects the interface $\Gamma_{fs}$. This implies that
$\T_h^a:=\{T\in\T_h\,:\, T\subset\Omega_a\}$, $a\in\{f,s\}$,  are admissible  triangulations of the fluid and poroelastic reference
domains $\Omega_a$, $a\in\{f,s\}$.
We exploit  the  finite element Taylor--Hood spaces which are popular in incompressible hydrodynamics:
\[
\begin{split}
\VV_h^a&=\{\bv\in C(\Omega_a)\,:\,\bv|_{T}\in P_2(T)^3\quad\forall\,T\in\T_h^a\},\quad a\in\{f,s\},\\
\QQ_h^a&=\{q\in C(\Omega_a)\,:\,q|_{T}\in P_1(T)\quad\forall\,T\in\T_h^a\},\quad a\in\{f,s\}.
\end{split}
\]
For trial functions we need also  the following subspaces:
\[
\begin{split}
\VV_h^{a,0}&=\{\bv\in\VV_h^a\,:\,\bv|_{\Gamma_{a,0}}=\mathbf{0}\},\quad a\in\{f,s\},\\
%\VV_h^{a,00}&=\{\bv\in\VV_h^{a,0}\,:\,\bv|_{\Gamma_{sf}}=\mathbf{0}\},\quad a\in\{f,s\},\\
\QQ_h^{s,0}&=\{q\in\QQ_h^s\,:\,q|_{\Gamma_{s,0}}={0}\}.\\
\end{split}
\]
We note that the Taylor--Hood is not a standard Darcy element for $H(div)$- formulations of the problem.
In particular, it fails to provide elementwise mass conservation. However, for applications where the local mass conservation is not a major concern,
it is a legitimate choice leading to optimal convergence in the Darcy region in product  $L^2$-velocity--$H^1$-pressure norm~\cite{karper2009unified}. %\textbf{Maxim, this is not correct, for the norms look at (2.4)-(2.5) in the paper}.
% coercivity of the Darcy-velocity  bilinear form on the subspace of weakly divergence free functions $Z_h^{a,0}=\{\bv\in\VV_h^{s,0}\,:\,\int_{\Omega_s}\Div\bv_h\,q_h\,dx={0}~\forall\,q_h\in\QQ^s_h\}$.
%To stabilize the formulation, we will add to a term penalizing the divergence of the Darcy velocity, an analogue of the grad-div stabilization for the Navier--Stokes problem.

For the time discretization, we assume a constant time step $\Delta t$ and use the notation $f^k(\bx)\approx f(k\Delta t, \bx)$ for all
time-dependent quantities.
 The first or second order  backward  finite difference approximation
 $\left[\frac{\partial f}{\partial t}\right]^{k}$
 of the time derivative of  $f$ at $t=k\Delta t$ is
\[
\left[\frac{\partial f}{\partial t}\right]^{k}=\frac{f^{k}-f^{k-1}}{\Delta t}\quad \text{or}\quad \left[\frac{\partial f}{\partial t}\right]^{k}=\frac{3 f^{k}-4 f^{k-1}+ f^{k-2}}{2\Delta t},
\]
respectively.
By $\widetilde{f}^k$ we denote the  extrapolated quantity  $f$
\begin{eqnarray*}
\widetilde{f}^k:=f^{k-1}\quad\text{or}\quad  \widetilde{f}^k:=2 f^{k-1}- f^{k-2}
\end{eqnarray*}
for the first or  second order extrapolation, respectively.
%We also use the short-hand notation $\bF_{k-\frac12}=\frac12(\bF_k+\bF_{k-1})$.

We proceed to multi-linear forms needed for our finite element formulation. For time derivatives, we need the form:
\begin{multline*}
 m^k(\bw_s,\bw_d,\bw_f,r;\bpsi_s,\bpsi_d,\bpsi_f,q):=
\int_{\Omega_s}\tilde J^{k}(\rho_p\bw_s+\rho_f\bw_d)\bpsi_s\,\mathrm{d}x
\\
+\int_{\Omega_s}\tilde J^{k}(\rho_f\bw_s+\frac{\rho_f}{\phi}\bw_d)\bpsi_d\,\mathrm{d}x
+
\int_{\Omega_s}s_0\tilde J^{k} r\,q\,\mathrm{d}x
+
\int_{\Omega_f}\rho_f\tilde J^{k}\bw_f\bpsi_f\,\mathrm{d}x.
\end{multline*}

For the elasticity part, we define
\[
a_s^k(\bw_s,\bpsi_s)=\int_{\Omega_s}\bF(\tbu^k)\bS(\bw_s,\tbu^k):\nabla\bpsi_s\,\mathrm{d}x\quad\text{and}\quad
a_d^k(\bw_d,\bpsi_s)=\int_{\Omega_s}\tilde J^{k} K^{-1}\bw_d\cdot\bpsi_d\mathrm{d}x,
\]
where $\bS(\bu_1,\bu_2)=\lambda_s\mbox{tr}(\bE(\bu_1,\bu_2))\bI+2\mu_{s}\bE(\bu_1,\bu_2)$,
$\bE(\bu_1,\bu_2)=\frac12\left\{\bF(\bu_1)^T\bF(\bu_2)-\bI\right\}_s$, $\{\bA\}_{s}=\frac12(\bA+\bA^T)$ denotes the symmetric part of
tensor $\bA\in\R^{3\times 3}$.

For the fluid domain we need the viscous term form
\[
a_f^k(\bw_f,\bpsi_f)=\int_{\Omega_f}2\mu_{f} \tilde J^k\bD_{\tbu^k}(\bw_f):\bD_{\tbu^k}(\bpsi_f)\,\mathrm{d}x
\]
and inertia form
\[
c_f^k(\bw_f;\bphi_f,\bpsi_f)=\int_{\Omega_f}\rho_f\tilde J^k\left(\nabla\bw_f\bF^{-1}(\tbu^k)
\bphi_f\right)\cdot \bpsi_f\,\mathrm{d}x,
\]
where $\bD_{\bu}(\bv)=\{(\nabla\bv)\bF^{-1}(\bu)\}_s$.

For handling the mass conservation constraints, we introduce
\[
b_a^k(q,\bpsi)=\int_{\Omega_a}q \tilde J^{k}\bF^{-T}(\tbu^k):\nabla \bpsi\,\mathrm{d}x, \quad a\in\{s,f\}.
%\quad\text{and}\quad b_f^k(q,\bpsi_f)=\int_{\Omega_f}q \tilde J^k \bF^{-T}(\tbu^k):\nabla \bpsi_f\,\mathrm{d}\bx.
\]

Next, we collect the interface  terms:
\[
\begin{split}
d^k(\bw_s,\bw_d,\bw_f,p_d;\bpsi_s,\bpsi_d,\bpsi_f)
=\tau\int_{\Gamma_{fs}} \tilde J^{k}_s \left((\bw_f-\bw_s-\bw_d)^T \bn\right)\left((\bpsi_f-\bpsi_s-\bpsi_d)^T \bn\right)\,ds \\
%d_{bjs}(\bw_s,\bw_f;\bpsi_s,\bpsi_f)
+ \int_{\Gamma_{fs}}\tilde J^{k}_s  p_d(\bpsi_f-\bpsi_s-\bpsi_d)\cdot \bn\,ds
%{\color{blue}- \int_{\Gamma_{fs}}\tilde J^{k}_s \frac{\rho_f}2(\bw_d\cdot \bn)_{-}\bw_f\cdot\bpsi_f\,ds}
+\int_{\Gamma_{fs}}\tilde J^{k}_s \frac{\rho_f}2|{\bv}_f|^2(\bpsi_s-\bpsi_f)\cdot\bn\,ds
 \\ +\gamma\int_{\Gamma_{fs}} \tilde J^{k}_s\ K^{-\frac12}\left(\bP(\bv_f-\bv_s)\right)\cdot\left(\bP(\bpsi_f-\bpsi_s)\right)\,ds ,
\end{split}
\]
with $\bn=\bF^{-T}\tbn/|\bF^{-T}\tbn|$, $\bP=\bI-\bn\bn^T$, and $ \tilde J^{k}_s= \tilde J^{k} |\bF^{-T}\tbn|$.
Parameter $\tau$ is a penalty parameter which forces the finite element solution to satisfy approximately  the normal velocity continuity condition.
The third   term on the right-hand side appears due to the additional term in the stress balance interface condition.

The finite element method with the backward difference time discretization reads:
Given $\bu^{k-1}$, $\bv^{k-1}_f$, $\bv^{k-1}_s$, $\bq^{k-1}$, $p^{k-1}_d$
find $\bv^{k}_f\in\VV_h^f$, $\bv^{k}_s\in\VV_h^{s,0}$, $\bq^{k}\in\VV_h^s$, $p^{k}_f\in\QQ_h^f$, $p^{k}_d\in\QQ_h^{s,0}$ such that
$\bv^{k}_f=\bg_h(\cdot,(k+1)\Delta t)$ on $\Gamma_{f0}$,
and the following identity holds:
\begin{equation}\label{FE1}
\begin{split}
 m^k\left({\small \left[\frac{\partial \bv_s}{\partial t}\right]^{k}},\right.&\left.{\small\left[\frac{\partial \bq}{\partial t}\right]^{k},\left[\frac{\partial \bv_f}{\partial t}\right]^{k},\left[\frac{\partial p_d}{\partial t}\right]^{k}};\bpsi_s,\bpsi_d,\bpsi_f,q_d\right)\\
 &+a_s^k(\bv^{k}_s,\bpsi_s)+a_d(\bq^{k},\bpsi_s)+a_f^k(\bv^{k}_f,\bpsi_f)+c_f^k(\bv^{k}_f,\widetilde{\bv}_f^k-\widetilde{\left[\frac{\partial \bu}{\partial t}\right]^{k}},\bpsi_f)\\
 &+d^k(\bv^{k}_s,\bq^{k},\bv^{k}_f;\bpsi_s,\bpsi_d,\bpsi_f) \\
 &-b_s^k(p_d,\bpsi_d)-b_f^k(p_f,\bpsi_f)+ b_s^k(q_d,\bv_s^k+\bq^k)+b_f^k(q_f,\bv_f^k)=0
 \end{split}
 \end{equation}
for all $\bpsi_f\in\VV_h^{f,0}$, $\bpsi_s\in \VV_h^{s,0}$,  $\bpsi_d\in \VV_h^{s}$,  $q_f\in\QQ_h^{f}$, $q_d\in\QQ_h^{s,0}$.
In addition, we relate the finite element displacement and the velocity field in the porous structure through the kinematic equation
\begin{equation}\label{FE2}
\left[\frac{\partial\bu}{\partial t}\right]^{k}=\bv^{k}_s\quad\text{in}~ \Omega_s.
\end{equation}
Equations \eqref{FE1}--\eqref{FE2} subject to the initial conditions and an equation for continuous extension of  $\bu^{k}$ from $\Omega_s$ onto
$\Omega_f$  define the discrete problem.
The continuous extension of $\bu$ in \eqref{ext} is provided by the elasticity equation   written for the velocity of the displacement~\cite{landajuela2017coupling}:

\begin{multline}\label{exten}
 -\Div\left[J\left(\lambda_m\text{tr}\left(\nabla\left[\frac{\partial\bu}{\partial t}\right]^{k}\mathbf{F}^{-1}\right)\mathbf{I}\right.\right. \\ \left.\left.+ \mu_m\left(\nabla\left[\frac{\partial\bu}{\partial t}\right]^{k}\mathbf{F}^{-1} + \left(\nabla\left[\frac{\partial\bu}{\partial t}\right]^{k}\mathbf{F}^{-1}\right)^{T}\right)\right)\mathbf{F}^{-T}\right] = 0~~~\text{in}~\Omega_f
 \end{multline}
satisfying the boundary condition $\left[\frac{\partial \bu}{\partial t}\right]^{k} = \bv^k$  on the interface $\Gamma_{fs}$.
The space dependent elasticity parameters are
$ \mu_m = \mu_s|\Delta_e|^{-1.2}$, $\lambda_m = 16\mu_m$,
where $|\Delta_e|$ denotes the physical volume of a mesh tetrahedron $\Delta_e$ subjected to displacement from the previous time step \cite{landajuela2017coupling}.

%\AV{Done with $\lambda_m$, $\mu_m$. The reference is changed!}.

Although the system is strongly coupled,  only \textit{a linear algebraic system should be solved on each time step}.

\section{Numerical experiments}\label{s_num}

In this section we assess the performance of the proposed monolithic FPSI FE method on the propagation of a pressure impulse in a compliant tube with
a porous wall filled with fluid. The problem setting follows
the benchmark suggested in \cite{pipe6} for flow in a tube with an impermeable  hyperelastic wall.
The original problem is related to the blood flow through an artery,
it has been extensively considered in the literature for validating the performance of FSI solvers~\cite{pipe3,pipe0,pipe4,pipe1,pipe2,pipe5}.
Since the test is an idealization of a practical setup, no experimental data is available and
the test serves to validate mesh convergence and study physical plausibility of the computed solution.

\begin{figure}[h!]
  \centering
 \begin{subfigure}{0.5\textwidth}
  \includegraphics[width=\textwidth]{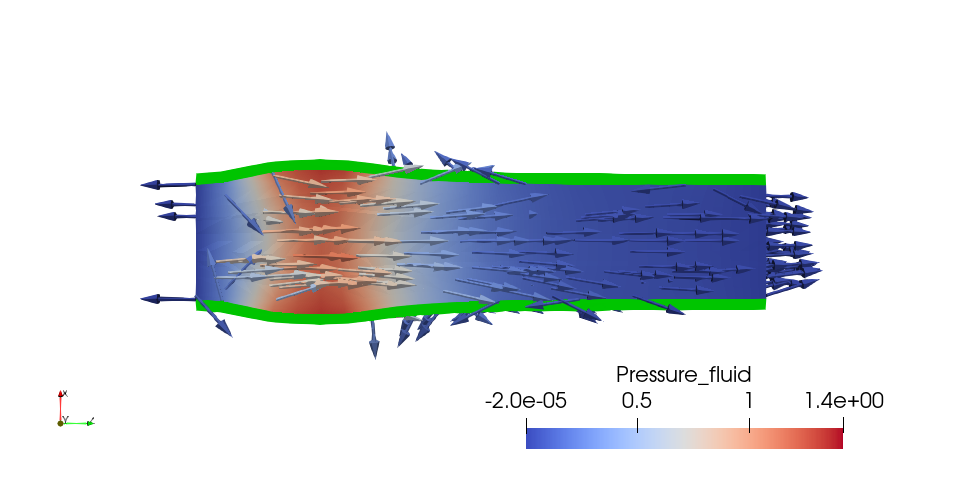}
  \caption{$t=0.004s$}
  \end{subfigure}\hskip0ex
   \begin{subfigure}{0.5\textwidth}
  \includegraphics[width=\textwidth]{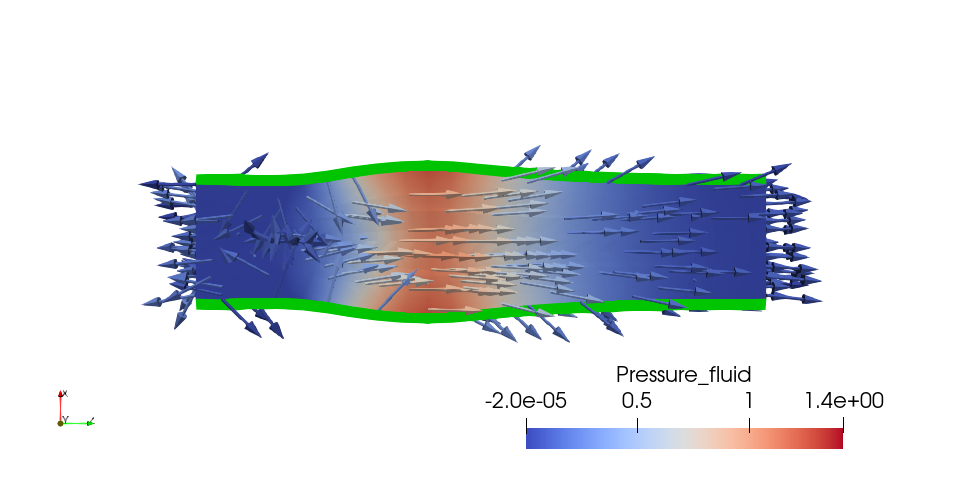}
  \caption{$t=0.006s$}
  \end{subfigure}
   \begin{subfigure}{0.5\textwidth}
  \includegraphics[width=\textwidth]{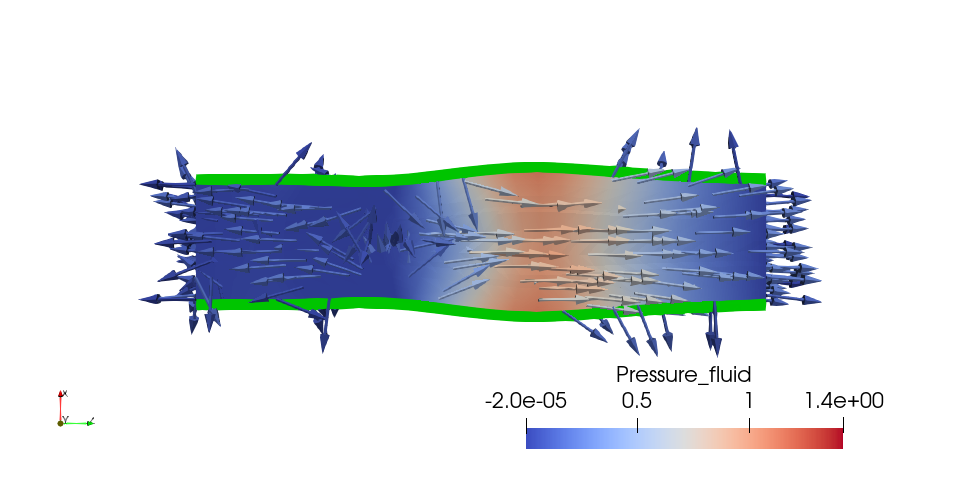}
  \caption{$t=0.008s$}
  \end{subfigure}\hskip0ex
  \begin{subfigure}{0.5\textwidth}
  \includegraphics[width=\textwidth]{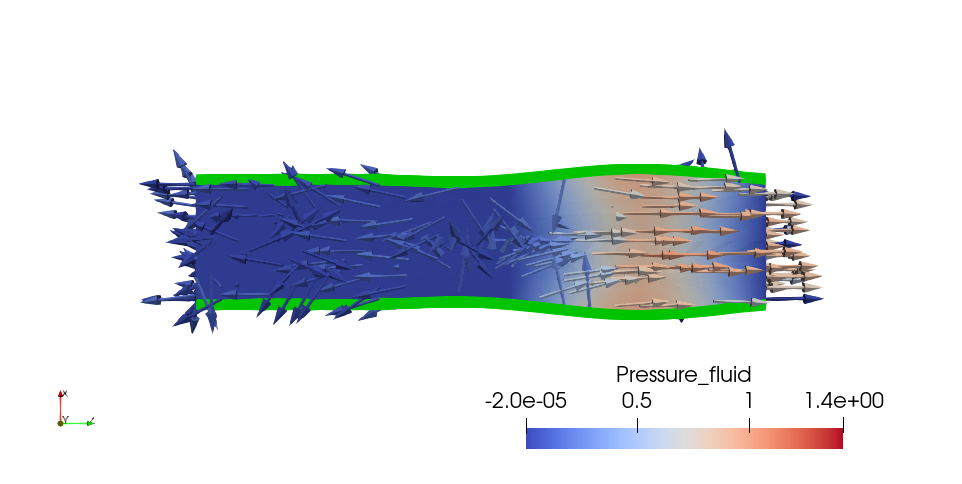}
  \caption{$t=0.01s$}
  \end{subfigure}
  \caption{Pressure wave: middle cross-section velocity field, pressure distribution, {velocity vectors} and 10-fold enlarged structure displacement  for several time instances.  \label{figure:Pipe1}}
\end{figure}

The problem configuration consists of an incompressible viscous flow through a poroelastic {tube with circular cross-section}. The tube is $50$mm long, it has  inner radius of $5$mm and the wall thickness is $1$mm.
The fluid density is {$10^{-3}$g/mm${}^3$} and kinematic viscosity is $3$mm${}^2$/s. The wall density  $\rho_s$ is $1.2\cdot10^{-3}$g/mm${}^3$.
In \eqref{constit_s1}, the Saint Venant--Kirchhof hyperelastic model is used with elastic modulus $E=3\cdot 10^{5}$g/mm/s${}^2$
and  Poisson's ratio  $\nu=0.3$. Initially, the fluid is at rest and the tube is non-deformed. The tube is fixed at both ends.

For the porous media parameters, we used porosity $\phi = 0.3$ \cite{rat},
mass storativity $s_0 = 5\cdot 10^{-5} mm\cdot s^2/g$ and two cases of the scalar
permeability coefficient: $K = 5\cdot10^{-13} mm^2$ and $K = 10^{-5} mm^2$.
The smaller value mimics permeability estimated in  rat's cardiovascular system \cite{rat}, while the larger value is taken from \cite{Yotov}.

On the left open boundary of the tube, the {external pressure $p_{ext}$} is set to {$1.333\cdot 10^{3}$}Pa for $t\in(0,3\cdot 10^{-3})$s and zero afterwards,
while on the right open boundary the {external pressure $p_{ext}$} is zero throughout the experiment.
This generates a pressure impulse that travels along the tube.
{The external pressure is incorporated into \eqref{FE1}--\eqref{FE2} through the open boundary condition $\bsigma_f\bF^{-T}\bn=p_{ext}\bn$.}

\begin{figure}[h!]
  \centering
 \begin{subfigure}{0.45\textwidth}
  \includegraphics[width=\textwidth]{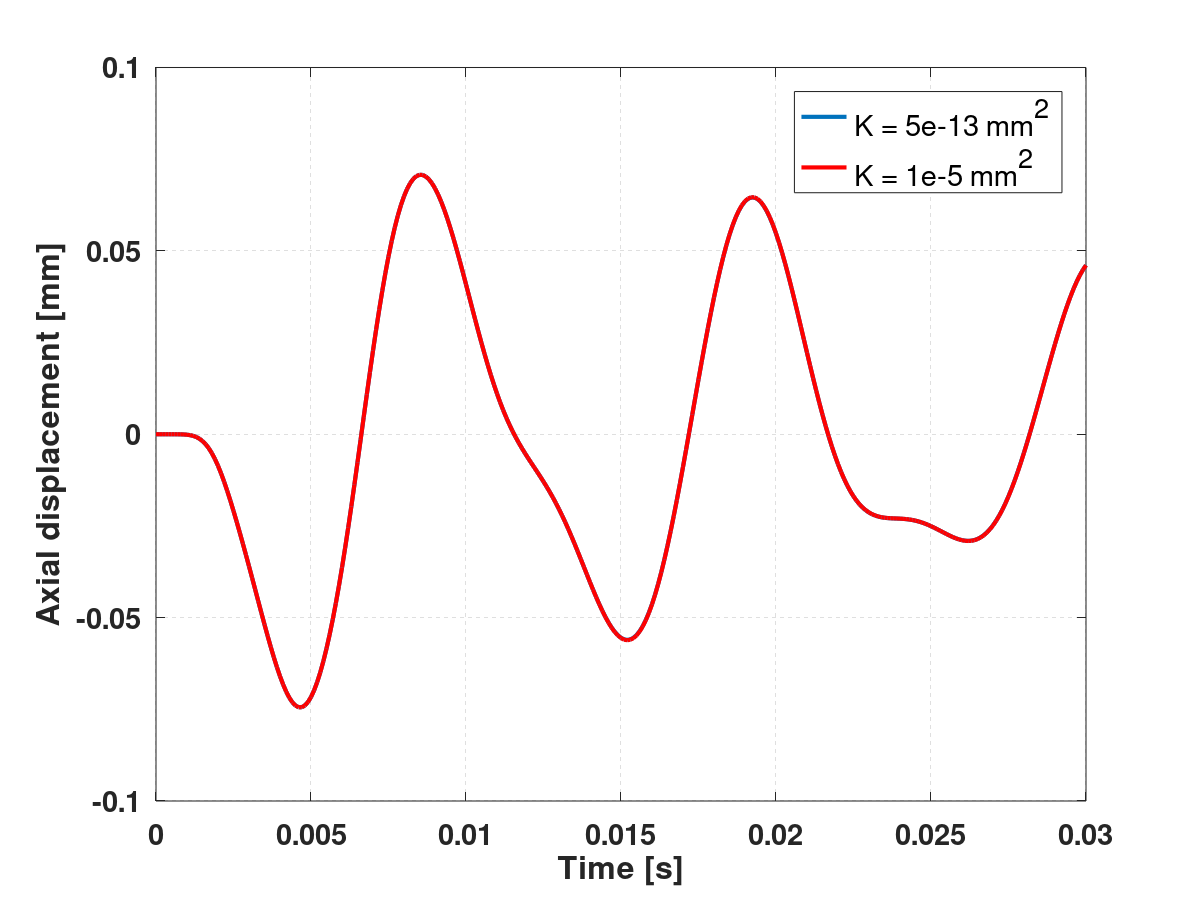}
  \caption{axial component}
  \end{subfigure}\hskip1ex
   \begin{subfigure}{0.45\textwidth}
  \includegraphics[width=\textwidth]{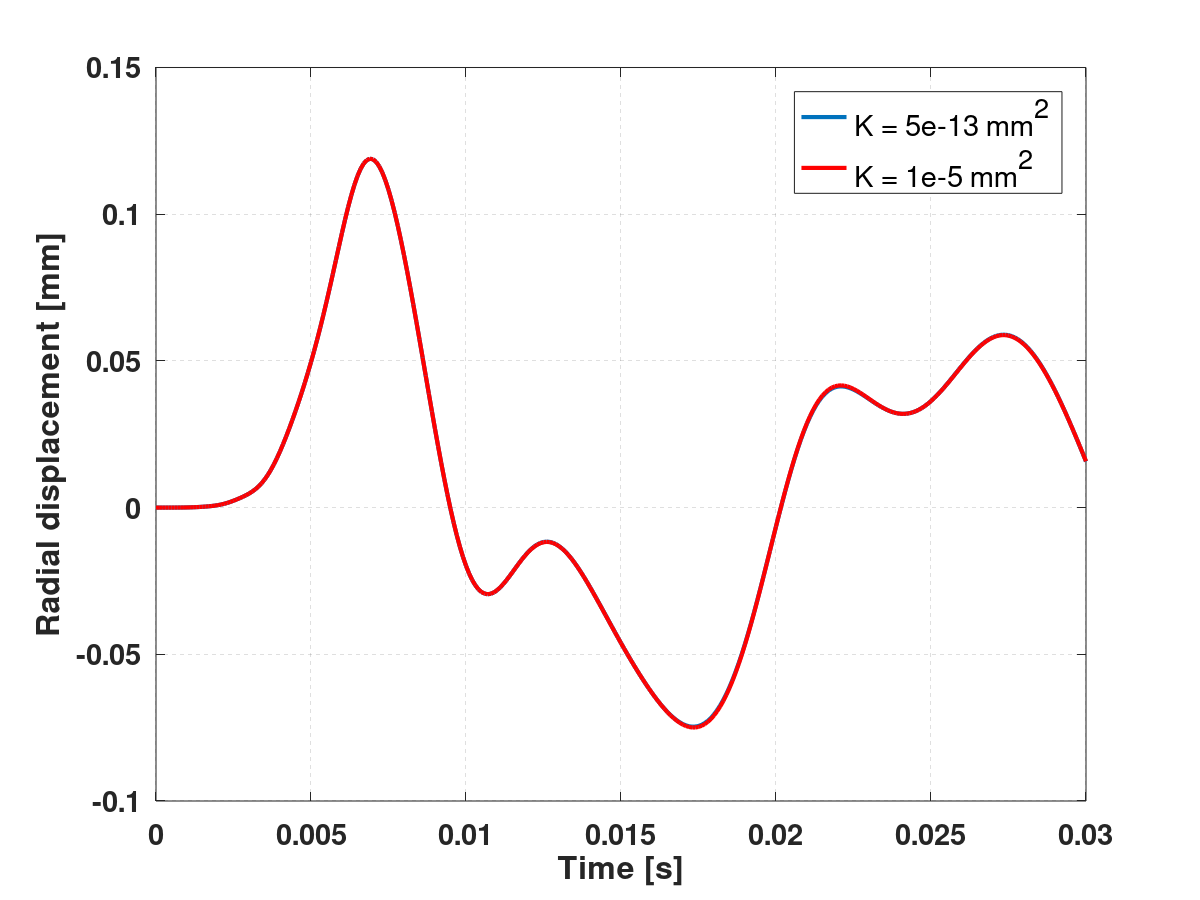}
  \caption{radial component}
  \end{subfigure}
  \caption{Pressure wave: The axial and radial components of displacement of the inner tube wall at half the length of the pipe. Solutions are shown for the two cases of permeability (see the text). {The plots are visibly indistinguishable.}  \label{figure:Pipe2}}
\end{figure}

\begin{figure}[h!]
  \centering
 \begin{subfigure}{0.5\textwidth}
  \includegraphics[width=\textwidth]{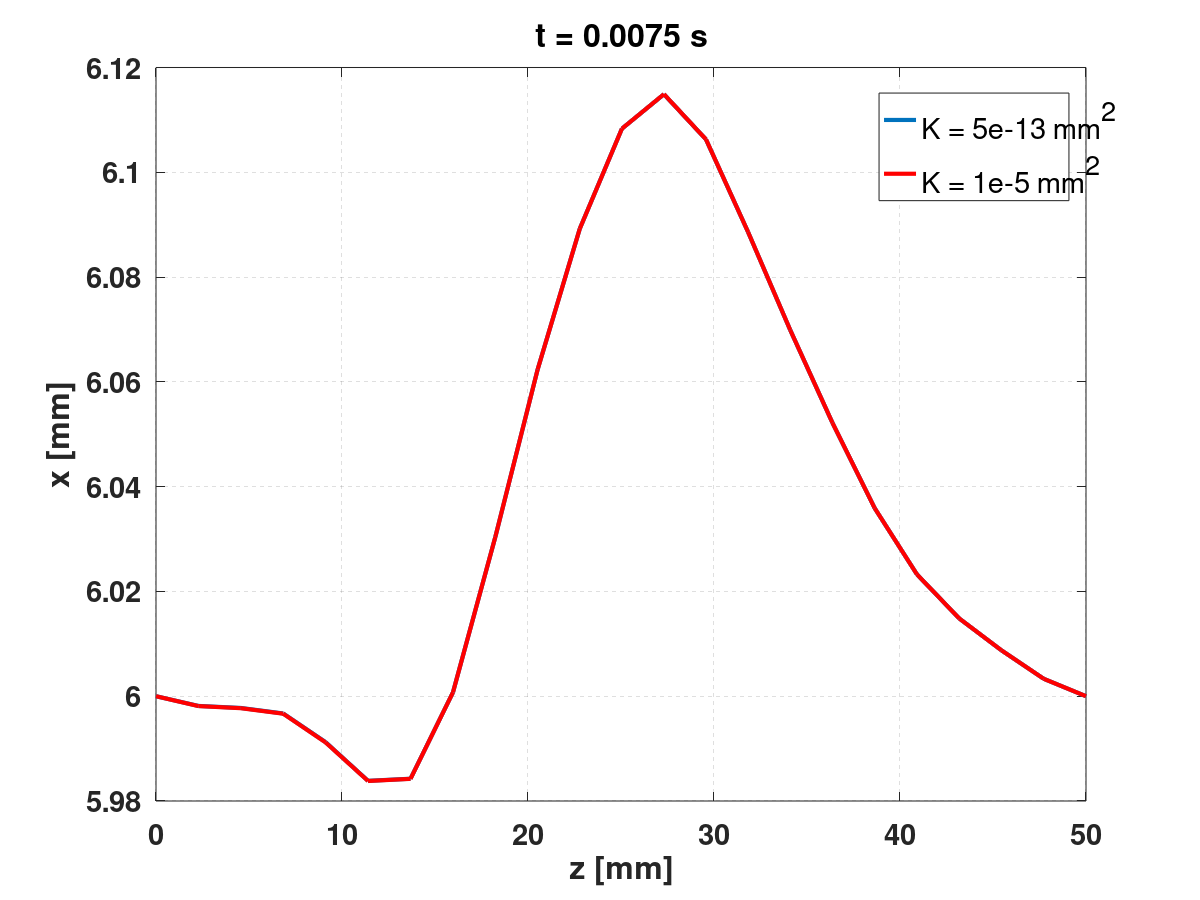}
  \caption{$t=0.004s$}
  \end{subfigure}\hskip0ex
   \begin{subfigure}{0.5\textwidth}
  \includegraphics[width=\textwidth]{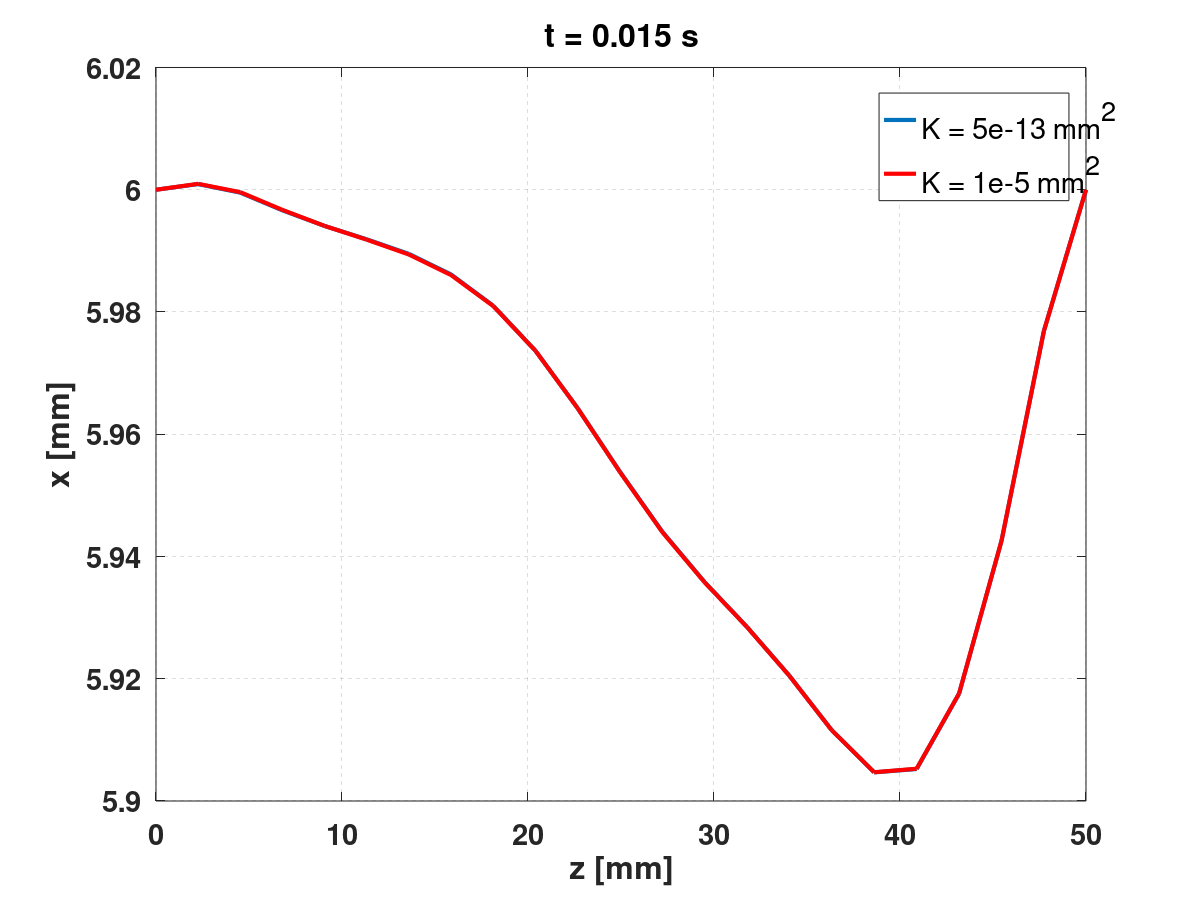}
  \caption{$t=0.006s$}
  \end{subfigure}
   \begin{subfigure}{0.5\textwidth}
  \includegraphics[width=\textwidth]{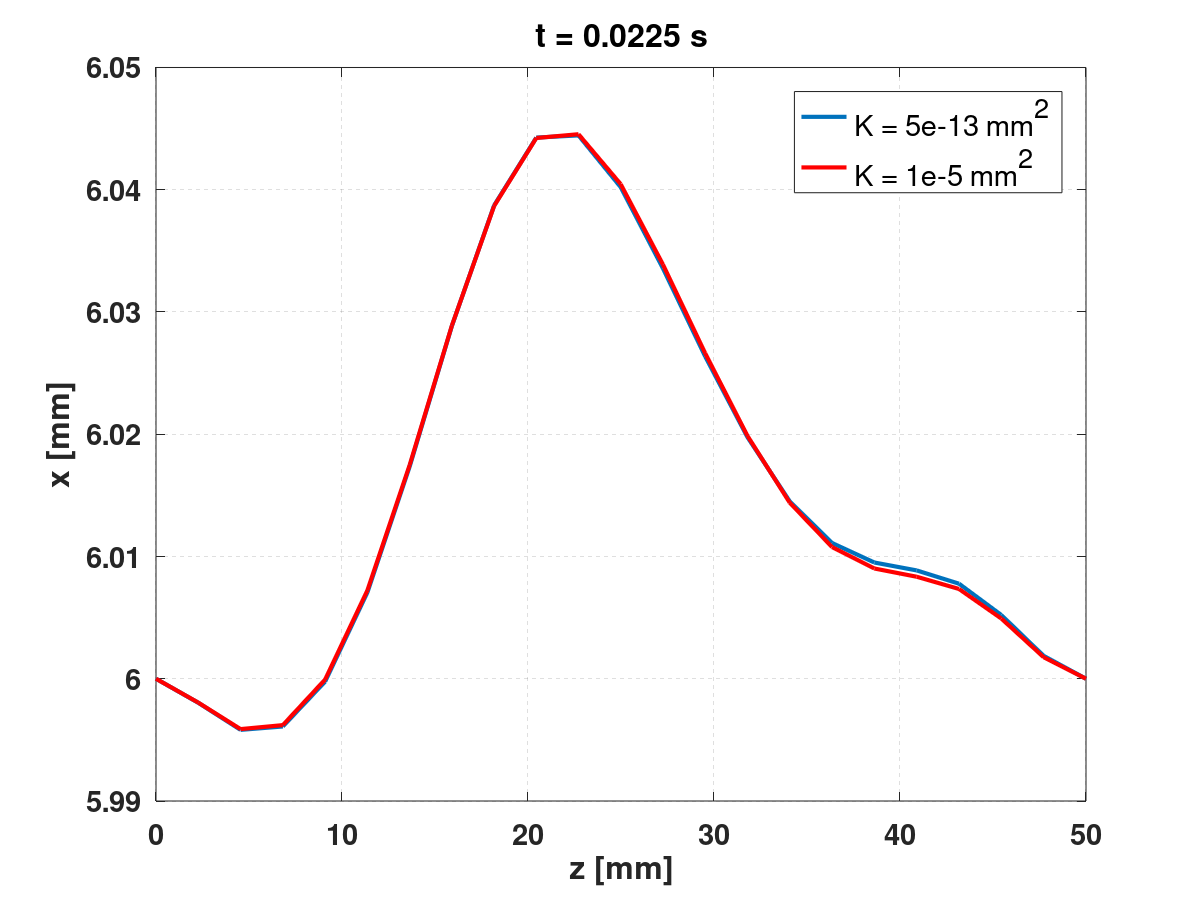}
  \caption{$t=0.008s$}
  \end{subfigure}\hskip0ex
  \begin{subfigure}{0.5\textwidth}
  \includegraphics[width=\textwidth]{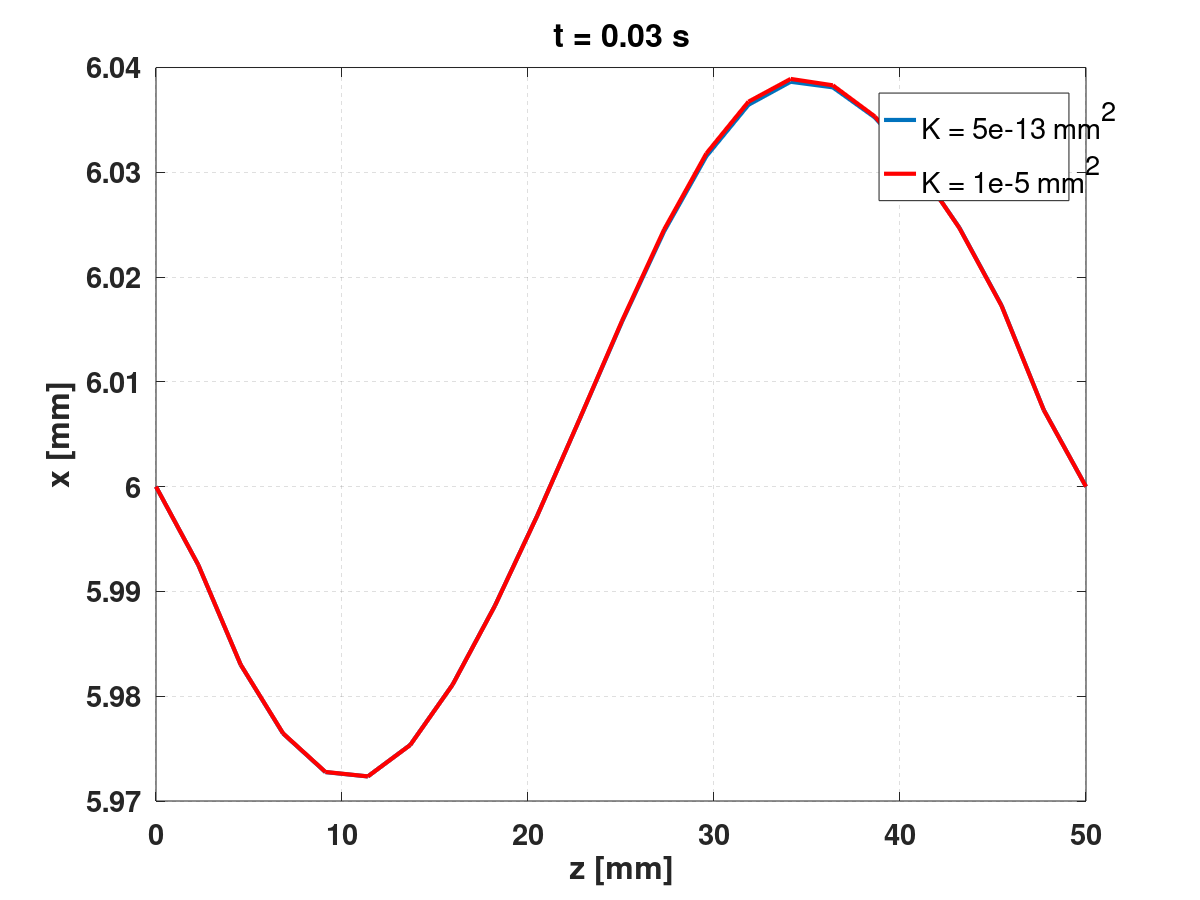}
  \caption{$t=0.01s$}
  \end{subfigure}
  \caption{Wall profile on the outer side along the tube length for several time instances.  \label{figure:profile}}
\end{figure}

%We built three computational meshes (coarse, fine and finer) for this experiment. On each level of refinement the mesh size was decreased by approximately a factor of $\sqrt{2}$.  The resulting numbers of tetrahedra in the fluid and solid subdomains,  respectively, are  13200 and 6336 for coarse, 29202 and 11904 for fine, and 89232 and 38016 for finer.

We use the Taylor--Hood P2-P1 elements for  velocity and pressure variables and P2 elements for displacements,
with the first order semi-implicit Euler discretization. The scheme \eqref{FE1}--\eqref{FE2} is implemented on the basis of
the open source package Ani3D \cite{Ani3D}.
The important feature of equation  \eqref{FE1} is linearization  on each time step due to extrapolation of all geometric factors and
the advection velocity  from the previous time steps. The resulting linear system  is solved by the multifrontal sparse direct solver MUMPS~\cite{MUMPS}.

The conformal mesh used for the numerical experiment has 13200 and 6336 tetrahedra for the fluid and solid subdomains, yielding  340586 degrees of freedom.
%, consisting of 14 blocks for the vertex degrees of freedom and 12 blocks for the edge degrees of freedom.
We set $\Delta t = 10^{-4}$s, $\gamma=1$, $\tau=h^{-2}$ where $h$ is the local mesh size. %For this problem, we do not use outflow boundary stabilization, i.e. $\kappa=0$ on both open boundaries.

%\AV{Instead, we incorporated boundary integrals of viscous stress tensor on both sides, although we did not observe noticeable difference compared to the computation without such integrals.} \MO{What does this exactly mean?}

Figure~\ref{figure:Pipe1} depicts  the computed fluid velocity field in the middle cross-section and wall displacement exaggerated by a factor of 10 for clarity. The redder the color of the arrow is, the larger magnitude the velocity vector has.

 %The results are similar to those found in other publications cited.
Figure~\ref{figure:Pipe2} shows the time variations of the radial and axial components of the displacement  of the inner tube wall
at half the length of the pipe, while Figure~\ref{figure:profile} shows the wall profile due to deformation at time instances $0.004, 0.006, 0.008, 0.01$.
Both Figures suggest that the difference in the permeabilities in this FPSI simulation scenario does not influence the FSI dynamics of the system.

%The maximum relative deviation between coarse and fine mesh displacement 2.1\% in the  axial component and  2.7\% in the radial component, which decreases to  maximum relative deviation between fine and finer mesh to 0.7\% in the axial component and  2.3\% in  the radial component.

%Both plots are consistent with the results reported in \cite{pipe6}.

\begin{figure}[h]
  \centering
 \begin{subfigure}{0.5\textwidth}
  \includegraphics[width=\textwidth]{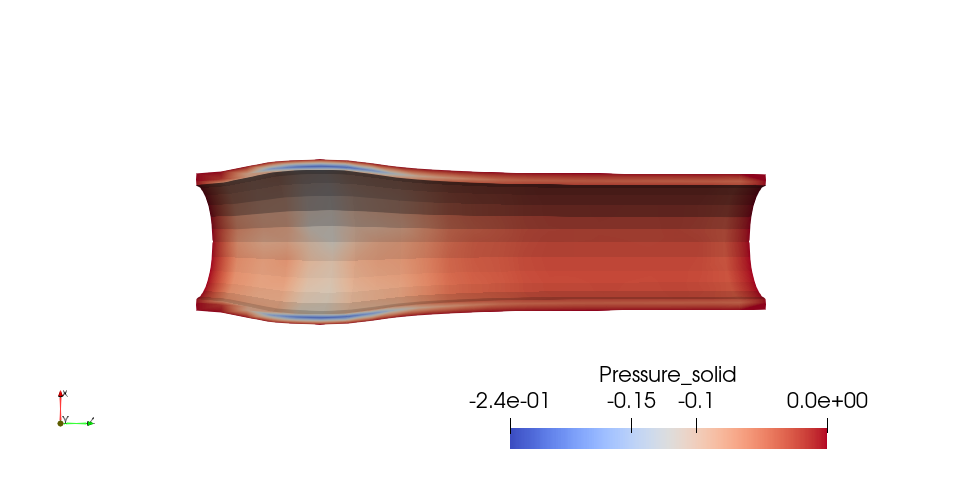}
  \caption{$t=0.004s$}
  \end{subfigure}\hskip0ex
   \begin{subfigure}{0.5\textwidth}
  \includegraphics[width=\textwidth]{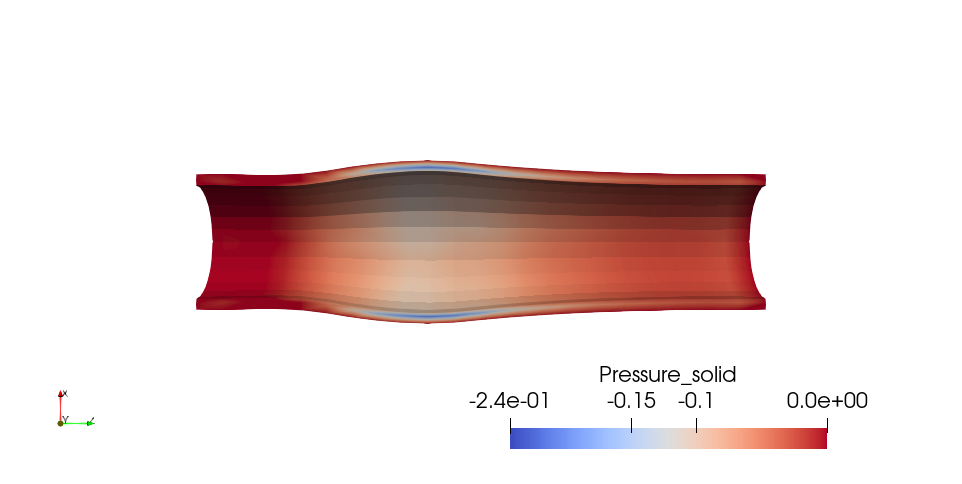}
  \caption{$t=0.006s$}
  \end{subfigure}
   \begin{subfigure}{0.5\textwidth}
  \includegraphics[width=\textwidth]{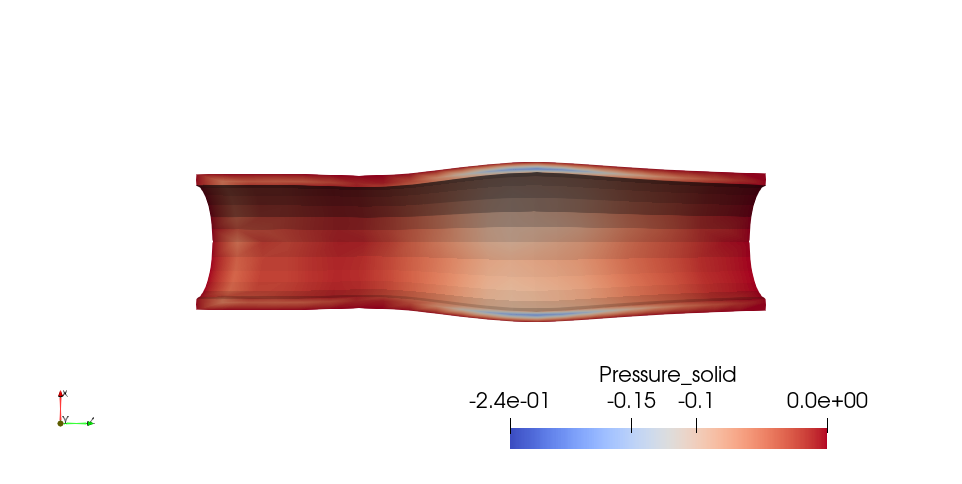}
  \caption{$t=0.008s$}
  \end{subfigure}\hskip0ex
  \begin{subfigure}{0.5\textwidth}
  \includegraphics[width=\textwidth]{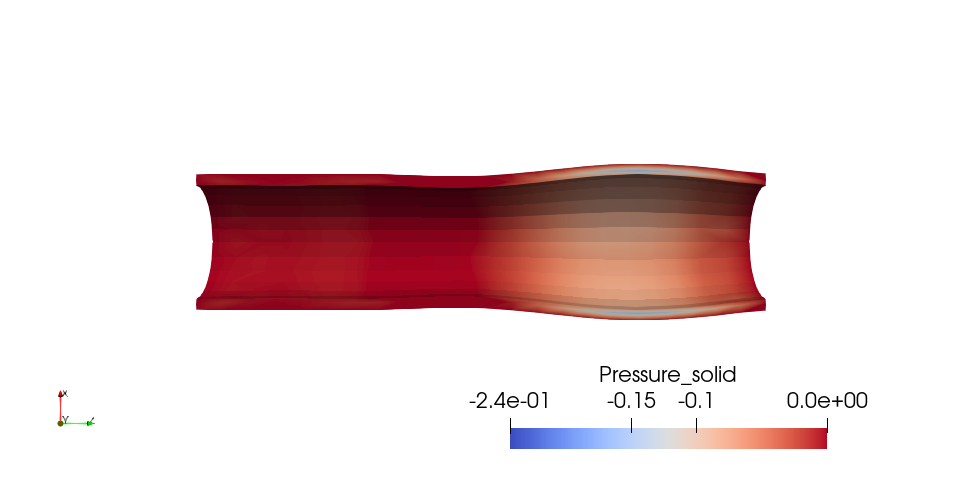}
  \caption{$t=0.01s$}
  \end{subfigure}
  \caption{Porous pressure $p_d$ distribution in the solid: middle cross-section view, with 10-fold enlarged structure displacement for several time instances.  \label{figure:Pipe3}}
\end{figure}

%\begin{figure}[h]
%  \centering
% \begin{subfigure}{0.4\textwidth}
%  \includegraphics[width=\textwidth]{q40.png}
%  \caption{$t=0.004s$}
%  \end{subfigure}\hskip2ex
%   \begin{subfigure}{0.4\textwidth}
%  \includegraphics[width=\textwidth]{q60.png}
%  \caption{$t=0.006s$}
%  \end{subfigure}
%   \begin{subfigure}{0.4\textwidth}
%  \includegraphics[width=\textwidth]{q80.png}
%  \caption{$t=0.008s$}
%  \end{subfigure}\hskip2ex
%  \begin{subfigure}{0.4\textwidth}
%  \includegraphics[width=\textwidth]{q100.png}
%  \caption{$t=0.01s$}
%  \end{subfigure}
%  \caption{Distribution of the magnitude of the porous velocity $\mathbf{q}$ in the solid: middle cross-section view, with 10-fold enlarged structure displacement for several time instances.  \label{figure:Pipe4}}
%\end{figure}

\begin{figure}[h]
  \centering
 \begin{subfigure}{0.5\textwidth}
  \includegraphics[width=\textwidth]{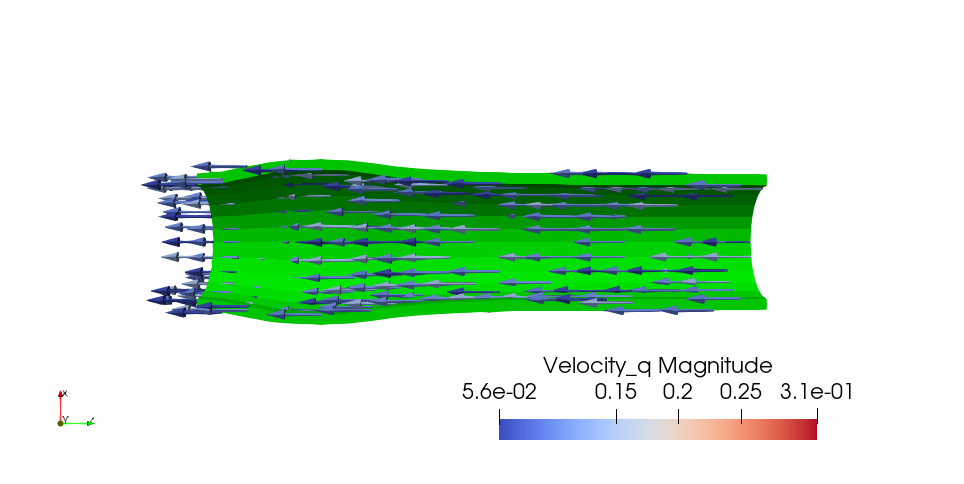}
  \caption{$t=0.004s$}
  \end{subfigure}\hskip0ex
   \begin{subfigure}{0.5\textwidth}
  \includegraphics[width=\textwidth]{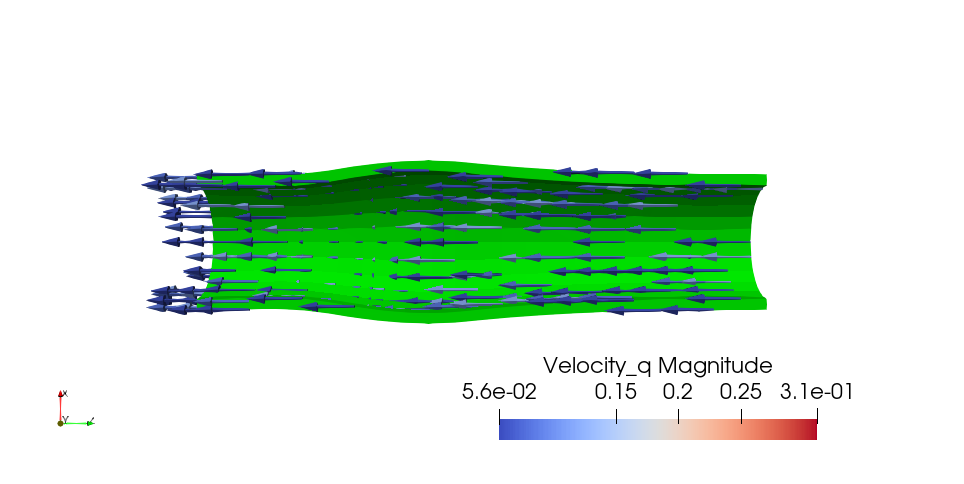}
  \caption{$t=0.006s$}
  \end{subfigure}
   \begin{subfigure}{0.5\textwidth}
  \includegraphics[width=\textwidth]{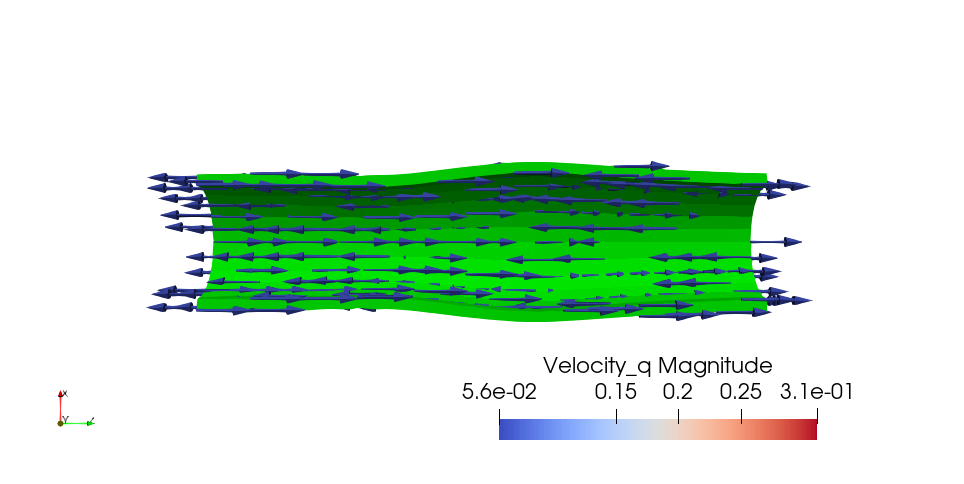}
  \caption{$t=0.008s$}
  \end{subfigure}\hskip0ex
  \begin{subfigure}{0.5\textwidth}
  \includegraphics[width=\textwidth]{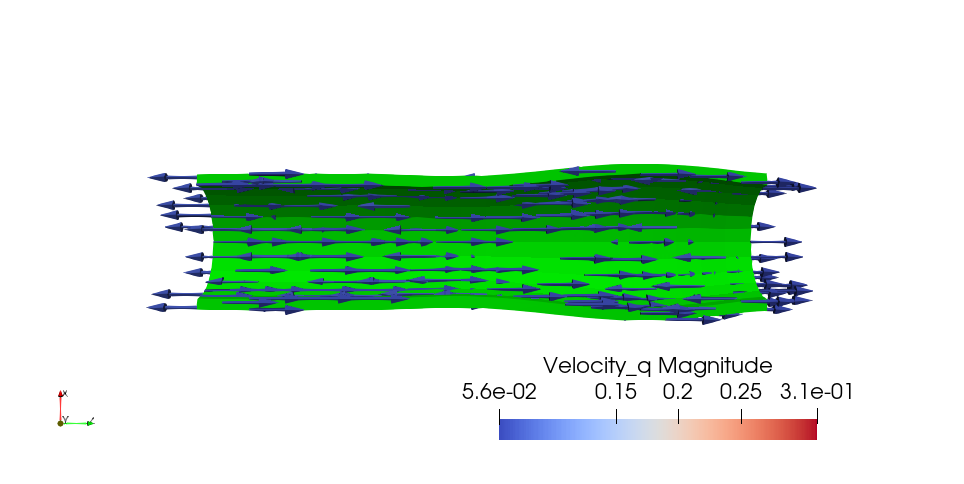}
  \caption{$t=0.01s$}
  \end{subfigure}
  \caption{Filtration velocity $\mathbf{q}$ distribution in the solid: middle cross-section view, with 10-fold enlarged structure displacement for several time instances.  \label{figure:Pipe5}}
\end{figure}

Figures~\ref{figure:Pipe3}-\ref{figure:Pipe5} demonstrate the porous pressure and filtration velocity distributions for
the same time instances $0.004, 0.006, 0.008, 0.01$ and permeability $K = 5\cdot10^{-13} mm^2$.
The maximum relative deviation for  filtration velocity $\mathbf{q}$ between the two permeability cases approaches 35\%:
the smaller permeability is, the  larger  magnitude of $\mathbf{q}$.
Both cases provide almost zero values for the non-axial components of $\mathbf{q}$.
The axial component of $\mathbf{q}$  points against the direction of the pressure pulse wave along the entire tube length.
The maximum relative deviation for the porous pressure $p_d$ is much lower, no more than 2.1\%.
The porous pressure is negative all across the tube wall and reaches zero value on the non-interface boundary according to the prescibed boundary conditions.

\clearpage

%\subsection*{Acknowledgements}
%The work was supported by the Russian Science Foundation grant 19-71-10094.

%\bibliographystyle{elsarticle-num}

% References
%\bibliographystyle{siam}
\bibliographystyle{plain}
\bibliography{mybib}
%\begin{thebibliography}{99}
%\end{thebibliography}

\end{document}